\documentstyle[amscd,amssymb,verbatim]{amsart}
\newcommand{\Tor}{\operatorname{Tor}}
\newcommand{\QGR}{\operatorname{QGR}}
\newcommand{\qgr}{\operatorname{qgr}}

\newcommand{\Vect}{\operatorname{Vect}}

\newcommand{\coker}{\operatorname{coker}}

\newcommand{\lan}{\langle}
\newcommand{\ran}{\rangle}

\newcommand{\CC}{{\cal C}}
\newcommand{\EE}{{\cal E}}

\newcommand{\Proj}{\operatorname{Proj}}

\newcommand{\si}{\sigma}

\newcommand{\Id}{\operatorname{Id}}

\renewcommand{\ker}{\operatorname{ker}}

\numberwithin{equation}{section}

\newtheorem{thm}{Theorem}[section]
\newtheorem{prop}[thm]{Proposition}
\newtheorem{lem}[thm]{Lemma}
\newtheorem{cor}[thm]{Corollary}
\newenvironment{defi}{\vspace{3mm}\noindent
{\bf Definition.}}{\vspace{3mm}}
\newenvironment{rem}{\vspace{3mm}\noindent
{\bf Remark.}}{\vspace{3mm}}
\newenvironment{rems}{\vspace{3mm}
\noindent {\bf Remarks.}}{\vspace{3mm}}

\newenvironment{exs}{\vspace{3mm}\noindent
{\bf Examples.}}{\vspace{3mm}}

\newcommand{\Pf}{\noindent {\it Proof}}

\newcommand{\ra}{\rightarrow}

\newcommand{\PP}{{\cal P}}

\newcommand{\OO}{{\cal O}}

\newcommand{\Hom}{\operatorname{Hom}}

\renewcommand{\a}{\alpha}

\newcommand{\Z}{{\Bbb Z}}

\newcommand{\Ga}{\Gamma}

\newcommand{\ot}{\otimes}

\newcommand{\sub}{\subset}
\newcommand{\ed}{\qed\vspace{3mm}}

\newcommand{\T}{\operatorname{Tor}} 
\newcommand{\coh}{\operatorname{coh}}
\newcommand{\cohproj}{\operatorname{cohproj}}
\newcommand{\GR}{\operatorname{GR}}

\title{Noncommutative $\Proj$ and coherent algebras}
\author{A.~Polishchuk}
\address
{Department of Mathematics,
University of Oregon, Eugene OR 97403} 
\email{apolish@@math.uoregon.edu}
\thanks{This work was partially supported by NSF grant DMS-0070967}

\begin{document}
\begin{abstract}
We prove that an abelian category equipped with an
ample sequence of objects is equivalent to the quotient of the
category of coherent modules over the corresponding algebra
by the subcategory of finite-dimensional modules. In the Noetherian
case a similar result was proved by Artin and Zhang in \cite{AZ}.
\end{abstract}

\maketitle

\centerline{\sc Introduction}

\bigskip

The main result of this paper is a slight generalization of the
theorem of Artin and Zhang in \cite{AZ} 
characterizing certain class of abelian categories that can be 
viewed as noncommutative analogues of categories of 
coherent sheaves on projective schemes. Recall that the main idea
of this approach to noncommutative projective geometry is to
associate to a noncommutative graded algebra $A$ the quotient
category $\QGR A$ of the category of graded $A$-modules  
by the subcategory of torsion modules. If $A$ is commutative and
is generated by a finite number of elements of degree $1$ then 
by the theorem of Serre \cite{S} the category of
quasicoherent sheaves on $\Proj(A)$ is equivalent to $\QGR A$.
Therefore, one would like to think about $\QGR A$ as a suitable
replacement for the latter category in the case when $A$ is
noncommutative. Theorem 4.5 of \cite{AZ} gives a nice criterion 
for a locally Noetherian abelian category $\CC$ to be equivalent to $\QGR A$
for some $A$. 
Namely, the criterion says that $\CC$ should
contain a sequence of objects $(E_n,n\in\Z)$ 
satisfying some properties analogous to those
of the sequence $(\OO_X(n),n\in\Z)$ for a projective scheme $X$. 
We will refer to this criterion
as the AZ-theorem and to $(E_n)$ as an {\it ample sequence}.

In the later developments of the 
above point of view on noncommutative geometry 
the assumption that the categories under consideration have to be 
(at least, locally) Noetherian was widely accepted as a convenient 
technicality. However, it appears
that the list of categories considered by ``noncommutative algebraic geometry''
should also contain some non-Noetherian categories. One family of such
examples is provided by the categories of holomorphic bundles on 
noncommutative tori. These categories can also  be viewed as hearts of
certain non-standard $t$-structures in
derived categories of coherent sheaves on elliptic
curves (see \cite{PS}). It is very easy to see that none of these categories
is Noetherian. More precisely, every non-zero object in these categories
is non-Noetherian. On the other hand, in the case when a noncommutative
two-torus
corresponds to a quadratic irrationality (has ``real multiplication'')
there is a natural graded algebra associated with this category and
one expects to have an analogue of Serre's theorem. 
The goal of the present paper is to prove that a (not
necessarily Noetherian) abelian category equipped with an ample sequence of 
objects can still be described in terms of the corresponding graded algebra.
Applications to noncommutative two-tori will be considered in \cite{P-rm}.

We exploit the idea going back to Serre's paper \cite{S} that the correct 
abelian category replacing the category of finitely generated modules in the
non-Noetherian case is the category of {\it coherent} modules.
Recall that an $A$-module $M$ is called coherent if $M$ is finitely
generated and the kernel of every surjection $A^{\oplus n}\ra M$ is
finitely generated. An algebra $A$ is called coherent if it is 
coherent as a module over itself. The graded version
of this definition can be reformulated as follows:
a finitely generated graded algebra $A$ is coherent iff the module of relations between
every finite collection of homogeneous elements in $A$ is finitely generated.
Roughly speaking, our main result is that if one removes the assumption that
the category is Noetherian in the AZ-theorem then the corresponding 
graded algebra is still coherent and the abelian
category in question is equivalent to the quotient of the category of
coherent modules by the subcategory of finite-dimensional modules. 

We did not study the question which coherent algebras appear in this way.
Recall that in the Noetherian case the answer is formulated
in terms of a rather unpleasant cohomological condition called $\chi_1$
(in the second part of the AZ-theorem).
It is not difficult to formulate a similar condition for coherent algebras
(replacing finitely generated modules by coherent modules in the definition
of \cite{AZ}) and we expect that an obvious extension of the second
part of the AZ-theorem is true in our situation.

One technical point: it is convenient to extend the class of graded algebras
to the wider class consisting of {\it $\Z$-algebras}. These are associative
algebras of the form $A=\oplus_{(i,j)\in\Z^2,i\le j}A_{ij}$ such that
the only non-zero products in $A$ are $A_{jk}\otimes A_{ij}\ra A_{ik}$.
The category of (nonnegatively) graded algebras is a
subcategory of the category of $\Z$-algebras: to every graded algebra
$A=\oplus_{n\ge 0} A_n$ one can associate a $\Z$-algebra $A_{\Z}=\oplus
A_{ij}$ with $A_{ij}=A_{j-i}$.
As was observed in \cite{SVdB}, sec. 11.1, the AZ-theorem
can be extended to the case when an ample sequence of objects does not 
have the form
$(\si^n(O),n\in\Z)$ for some object $O$ and some autoequivalence $\si$, by working
with $\Z$-algebras. Our generalization is also formulated using the language of
$\Z$-algebras. 

It is worth mentioning that the notion of a coherent
$\Z$-algebra arises naturally in the theory of geometric helices
developed in \cite{BP}.
Namely, it was proved in \cite{P} that in the situation when a
triangulated category is generated by the geometric helix,
certain natural pair of 
subcategories defines a $t$-structure iff the $\Z$-algebra
associated to this helix is coherent. 

It would be very desirable to develop some other techniques for
checking whether a given graded algebra is coherent.
We present some partial results in this direction which allow
us to construct various examples of coherent and noncoherent algebras. 
Also, as a consequence of our main result we derive the connection
between coherency of a graded algebra and its Veronese subalgebras.
We believe that in the noncommutative world many natural construction lead to
coherent algebras (but not necessarily Noetherian ones). For example,
Piontkovskii in \cite{Pi1} proved coherence of a graded algebra with finite number of
generators and a finite number of defining {\it monomial} relations. In \cite{Pi2}
this result is generalized to a broader class of algebras.

{\it Acknowledgments}. This note is mostly based on a part of the author's diploma work \cite{P} 
carried out at Moscow State University in 1993 under the supervision of A.~Bondal
to whom I am very grateful. I also thank L.~Positselski and
D.~Piontkovskii for helpful discussions. 

{\it Conventions}. Throughout this paper we work over a fixed field $k$. By a {\it graded algebra}
we mean a nonnegatively graded $k$-algebra of the form 
$A=\oplus_{i\ge 0}A_i$, where $A_0=k$ (thus, we consider only {\it connected}
algebras).

\section{Preliminaries on $\Z$-algebras}

In this section we review some basic constructions dealing with
$\Z$-algebras.

\begin{defi}(see \cite{BGS},\cite{BP}) A {\it $\Z$-algebra} is
an associative $k$-algebra $A$ of the form 
$A=\oplus_{i\le j}{A_{ij}}$, $i,j\in \Z$, where
$A_{ii}=k$ for all $i\in\Z$; the only non-zero components of multiplication with
respect to this decomposition are $A_{jk}\ot A_{ij}\ra A_{ik}$;
moreover, the multiplications $A_{jj}\ot A_{ij}\ra A_{ij}$
and $A_{ij}\otimes A_{ii}\ra A_{ij}$ are the identity maps.
\end{defi}

We will always impose the following finiteness condition on
a $\Z$-algebra: $\dim_k A_{ij}<\infty$ for all $i\le j$.
We consider the category $\GR A$ of graded right $A$-modules
$M=\oplus_{i\in\Z}{M_i}$ with an $A$-action of the form 
$M_j\ot A_{ij}\ra M_i$, such that $M_i\ot A_{ii}\ra M_i$
are the identity maps. The morphisms in $\GR A$ are homomorphisms
of $A$-modules (they preserve the grading automatically).
We denote by $S_j$ the unique irreducible $A$-module such that
$(S_j)_i=0$ for $i\neq j$, $(S_j)_j=k$ and by $P_j=\oplus_i{A_{ij}}$
its indecomposable projective cover. Let $\PP$ be the family of modules
consisting of the finite direct sums of $P_j$'s. We say that
an $A$-module $M$ is {\it finitely generated} 
(resp., {\it finitely presented}) if
there is a surjection $P\ra M$ with $P\in \PP$ (resp., 
$M=\coker(P'\ra P)$, where $P,P'\in \PP$). Note that our assumption on
$A$ implies that for every finitely generated $A$-module 
$M=\oplus_i M_i$ the graded components $M_i$ are finite-dimensional.

As we already mentioned before, to a graded algebra 
$A=\oplus_{i\ge 0}{A_i}$ with $\dim_k A_i<\infty$ one can
associate the $\Z$-algebra $A_{\Z}=\oplus{A_{ij}}$ with
$A_{ij}=A_{j-i}$. Note that the categories of graded right modules
over $A$ and $A_{\Z}$ are equivalent: the equivalence sends
an $A$-module $M=\oplus M_i$ to an $A_{\Z}$-module $M_{\Z}=\oplus
M_{-i}$. Under this equivalence $P_j$ corresponds to the 
free module $A(j)$,
where $A(j)_i=A_{i+j}$.

Let $A$ be a $\Z$-algebra. Henceforward, by an $A$-module we 
always mean a graded right $A$-module. 

\begin{defi} An  $A$-module $M$ is called {\it coherent} if it 
satisfies the following two conditions: 

\noindent
(i) $M$ is finitely generated; 

\noindent
(ii) for every homomorphism $f: P\ra M$ with $P\in \PP$ the
module $\ker(f)$ is finitely generated.
\end{defi}

Of course, this definition is essentially a particular case of
the general definition of coherent sheaves of modules given
by Serre in \cite{S}.
We denote by $\coh A\subset \GR A$ the full subcategory of
coherent modules. Some basic properties of coherent modules
proved in \cite{S} also hold in our situation with same proofs. 
Most notably, we have the following result.

\begin{prop}\label{sprop} $\coh A$ is an abelian
subcategory of $\GR A$ closed under extensions.  
\end{prop}

\begin{defi} A $\Z$-algebra $A$ is called {\it weakly right coherent} if all
the modules $P_j$ are coherent. It is called 
{\it right coherent} if in addition
all the modules $S_j$ are coherent.
\end{defi}

Similarly, one defines the notion of a left coherent $\Z$-algebra.
In this section and in the next one
we work exclusively with right modules, so
by coherence we mean right coherence. 
Considering graded algebras as a subcategory in the category of
$\Z$-algebras we get a notion of coherence for them. 
Note that for a finitely generated graded algebra weak coherence and coherence
are equivalent. Later we will prove
that the tensor algebra $T(V)$ of a finite-dimensional space $V$
is coherent (see Corollary \ref{tensorcor}). This easily implies 
that the tensor algebra $T(V)$ of an infinite-dimensional
vector space is weakly coherent (but not coherent).

\begin{lem}\label{reslem} Assume that $A$ is weakly coherent.
Then every coherent module $M$ has a resolution 
$\ldots P^{-2}\ra P^{-1}\ra P^0\ra M$
with $P^i\in\PP$.
\end{lem}

\Pf . Indeed, since $M$ is finitely generated, 
we can choose a surjection $f:P^0\ra M$,
where $P^0\in\PP$. By Proposition \ref{sprop}
$\ker(f)$ is again coherent, so we can iterate
this procedure.
\ed

One can easily prove the following criterion: a $\Z$-algebra $A$ is weakly
coherent (resp., coherent) 
iff there exists a full abelian subcategory $\CC\subset \GR A$ 
consisting of finitely generated modules, such that $\PP\subset
\CC$ (resp., $\PP\sub\CC$ and $S_j\in\CC$ for every $j$). 
This criterion can be used to prove the following result.

\begin{prop}\label{homcohprop} 
Let $A\ra B$ be a homomorphism of $\Z$-algebras and
$\{ Q_j,\ j\in \Z \}$ be the set of
indecomposable projective $B$-modules (defined in the same way as
$P_j$ for $A$). Assume that $A$ is coherent and all $Q_j$ are
coherent as $A$-modules. Then $B$ is a coherent $\Z$-algebra.
\end{prop}

\Pf . This follows immediately from the above criterion 
applied to the subcategory
of $B$-modules that are coherent over $A$. 
\ed

For example, the above proposition implies
that the quotient $A/J$ of a coherent graded algebra
$A$ by a two-sided (homogeneous) ideal $J$, such that $J$ is 
finitely generated as a right ideal, is again a coherent algebra.

A module $M=\oplus_i M_i$ 
over a $\Z$-algebra $A$ is called {\it bounded} (resp., {\it bounded
above}) if $M_i\neq 0$ for only a finite number of indices $i$
(resp., $M_i=0$ for all sufficiently large  $i$).
Note that every finitely generated $A$-module is bounded above
(since this is true for $P_j$).
For an $A$-module $M=\oplus_i M_i$ and an integer $n\in\Z$ 
we denote by $M_{\le n}$ the $A$-module
$\oplus_{i\le n} M_i$.
For a coherent $\Z$-algebra $A$ we denote by 
$\coh^b A$ the category of bounded coherent
modules. Note that since all the modules $S_j$ are coherent,
$\coh^b A$ consists exactly of all finite-dimensional modules.
Furthermore, $\coh^b A$ is a Serre subcategory of $\coh A$, so
the quotient category 
$$\cohproj A:=\coh A/ \coh^b(A)$$ 
is still abelian. In particular, we can construct such
a category (still denoted $\cohproj A$) for a coherent
graded algebra $A$. If
$A$ is Noetherian then $\cohproj A$ coincides with the category $\qgr A$
considered in \cite{AZ}.

\section{Coherent sequences and an equivalence of categories}

Let $\CC$ be an abelian $k$-linear category, 
$\EE=(E_i,\ i\in \Z)$ be a sequence of objects of $\CC$. 
We are going to show that under appropriate conditions on $\EE$
(that are combined below in the notion of a {\it coherent sequence})
certain quotient of the category $\CC$ is equivalent to
$\cohproj A$, where the  
$\Z$-algebra $A=A(\EE)=\oplus_{i\le j} {A_{ij}}$ is defined as follows:
$A_{ij}=\Hom_{\CC}(E_i,E_j)$ for $i<j$, $A_{ii}=k$, the multiplication
is induced by the composition in $\CC$.
In the particular case of an {\it ample} sequence, we will obtain an
equivalence of $\CC$ with $\cohproj A$. 

Using this technique we will show that for a coherent
graded algebra $A$ generated by $A_1$ over $A_0=k$ there is an
equivalence of categories $\cohproj A\simeq \cohproj A^{(n)}$,
where $A^{(n)}=\oplus_{i\ge 0}{A_{in}}$ is a Veronese subalgebra of
$A$. The similar result in the Noetherian
case is due to Verevkin \cite{Ver}, Theorem (A-5) 
(see also \cite{AZ}, Prop. 5.10).

From now on we will always assume that the sequence $\EE=(E_i)$ satisfies
the following finiteness condition:
for every object $X\in\CC$ one has $\dim_k\Hom_{\CC}(E_i, X)<\infty$.
In particular, the components $A_{ij}$ of the corresponding $\Z$-algebra
$A$ are finite-dimensional.

\begin{defi} A sequence $\EE=(E_i)$ is called {\it projective} if for 
every surjection $X\ra Y$ in $\CC$ there exists $n\in\Z$
such that the corresponding map
$\Hom(E_i,X)\ra \Hom (E_i,Y)$ is surjective for $i<n$.
\end{defi}

Every sequence $\EE=(E_i)$ defines a functor 
$\Ga_*:\CC\ra \GR A$ sending $X$ to the right $A$-module 
$\oplus_{i\in\Z}\Hom(E_i,X)$. 
More important for us are truncated versions of this functor
$$\Ga_{\le m} X:=(\Ga_* X)_{\le m}=\oplus_{i\le m}\Hom(E_i,X),$$
where $m\in\Z$. 
If $\EE$ is projective then $\Ga_*$ (resp., $\Ga_{\le m}$) is
exact modulo the subcategory of bounded below (resp., bounded) modules. 
We define the subcategory $\CC_0=\CC_0(\EE)\subset
\CC$ as the full subcategory consisting of objects $X$ such that
$\Hom(E_i,X)=0$ for $i<<0$. If $\EE$ is projective then $\CC_0$ is a Serre
subcategory.

\begin{defi} (i) A projective sequence $\EE=(E_i)$ is called {\it coherent} if for
every object $X\in\CC$ and every $m\in\Z$ there exists a set of integers
$i_1, \ldots, i_s$ with $i_j\le m$ for all $j$, such that the canonical map
$$\oplus_{j=1}^s\Hom(E_{i_j},X)\ot \Hom(E_i, E_{i_j})\ra \Hom(E_i,X)$$
is surjective for $i<<0$.

\noindent (ii) A coherent sequence $\EE=(E_i)$ is called {\it ample} if $\CC_0(\EE)=0$.
\end{defi}

The following lemma gives convenient reformulations of these conditions. 

\begin{lem}\label{cohdeflem} 
Let $\EE=(E_i)$ be a projective sequence.

\noindent
(i) $\EE$ is coherent iff for every object $X\in\CC$ the $A$-modules 
$\Ga_{\le m}X$ are finitely generated for all $m$.

\noindent
(ii) $\EE$ is ample iff for every $X\in \CC$ and every $m\in\Z$ 
there exists a surjection $\oplus_{j=1}^s E_{i_j}\ra X$ for some
$i_1,\ldots,i_s$ with $i_j\le m$ for all $j$.
\end{lem}

\Pf . (i) Assume first that $\EE$ is coherent. Note that
for all $i\le m$ the module $\Ga_{\le m}E_i$ contains $P_i$ as a submodule and
the quotient is finite-dimensional, hence, $\Ga_{\le m}E_i$ is finitely
generated in this case. Now for an object $X\in \CC$ and an
integer $m$ we choose integers $i_1,\ldots, i_s\le m$ as in the definition
and consider the morphism 
\begin{equation}\label{cohmorphism}
f:\oplus_{j}\Hom(E_{i_j},X)\ot E_{i_j}\ra X,
\end{equation}
such that $\coker \Ga_*f$ is bounded below. Since by our assumption on $\EE$
all spaces $\Hom(E_i,X)$ are finite-dimensional, this implies that 
$\coker \Ga_{\le m}f$ is finite-dimensional. On the other hand, 
as we observed above, the $A$-modules $\Ga_{\le m}E_{i_j}$ are finitely generated. 
It follows that $\Ga_{\le m}X$ is finitely generated.
Conversely, if $\Ga_{\le m}X$ is finitely generated then there exists a surjection
$\oplus_{j=1}^s P_{i_j}\ra\Ga_{\le m}X$, where $i_j\le m$ for all $j$, hence,
the condition in the definition of coherence is satisfied for $m$.

\noindent
(ii) Assume that $\EE$ is ample. Then for every $m$ we can
choose integers $i_1,\ldots,i_s\le m$
as in the definition of coherence and consider the morphism (\ref{cohmorphism}).
We claim that $\coker f=0$. Indeed, since $\EE$ is projective, we have 
$\Ga_{\le n}(\coker f)\simeq\coker\Ga_{\le n} f=0$ for $n<<0$. Hence, $\coker f$
belongs to $\CC_0=0$. Conversely, assume that for every $X\in\CC$ and $m\in\Z$ 
there exists a surjection $f:\oplus_{j=1}^s E_{i_j}\ra X$, where $i_j\le m$ for
all $j$. Then by projectivity of $\EE$ the map $\Ga_{\le n}f$ is surjective for $n<<0$
which implies that the sequence $\EE$ is coherent. On the other hand, it is clear that
in this case $\CC_0=0$.
\ed

\begin{cor}
If $A$ is a coherent $\Z$-algebra then the sequence $(P_i)$ in $\coh A$ is coherent.
\end{cor}

\Pf . Indeed, all truncations of a coherent module are still coherent.
\ed

\begin{rem} The above lemma shows that our definition of ampleness is a direct generalization
of the definition (4.2.1) of \cite{AZ}.
\end{rem}

\begin{prop}\label{cohprop} Let $\EE=(E_i)$ be a coherent sequence, $A=A(\EE)$ be
the corresponding $\Z$-algebra. Then

\noindent
(i) for every $X\in \CC$ and every $m\in\Z$ the $A$-module $\Ga_{\le m}X$ is coherent;

\noindent
(ii) the $\Z$-algebra $A$ is coherent.
\end{prop}

\Pf . (i) From Lemma \ref{cohdeflem}(i) we know that the modules $\Ga_{\le m}X$ are coherent.
Let $f:\oplus_{j=1}^s P_{i_j}\ra \Ga_{\le m}X$ be any
homomorpism of $A$-modules. We have to show that $\ker f$ is finitely generated. 
Clearly we can assume that $i_j\le m$ for all $j$. Then $f$
corresponds to a morphism $\phi:\oplus_j E_{i_j}\ra X$ in $\CC$. Let $n$ be an integer
smaller than all $i_j$. We have
$$(\ker f)_{\le n}=\ker\Ga_{\le n}\phi\simeq\Ga_{\le n}(\ker \phi),$$ 
so by Lemma \ref{cohdeflem}, this module is finitely generated.
Hence, $\ker f$ is also finitely generated.

\noindent 
(ii) By part (i) the modules 
$(P_j)_{\le j-1}=\Ga_{\le j-1}E_j$ and $\Ga_{\le j}E_j$
are coherent for every $j$. Therefore, the
exact triple
$$0\ra \Ga_{\le j-1}E_j\ra \Ga_{\le j}E_j\ra \Hom(E_j,E_j)\ot S_j\ra 0$$
implies that $S_j$ is a coherent $A$-module. 
But $P_j$ is an extension of
$S_j$ by $(P_j)_{\le j-1}$, so it is also coherent. 
\ed

Now we can prove our main theorem that generalizes the first part of Corollary 4.6
in \cite{AZ}.

\begin{thm}\label{mainthm} 
Let $\EE=(E_i)$ be a coherent sequence, $A=A(\EE)$ be the
corresponding algebra. Then there is an equivalence of categories
$\CC/\CC_0 \simeq \cohproj A$. In particular, if $\EE$ is ample
then $\CC\simeq \cohproj A$.
\end{thm}

\Pf . Recall that according to Proposition
\ref{cohprop} the $\Z$-algebra $A$ is coherent and for every $X\in \CC$
and every $m\in\Z$ the $A$-module $\Ga_{\le m} X$ is coherent.
Now any of the functors $\Ga_{\le m}:\CC\ra \coh A$ 
induces an exact functor 
$$\Phi:\CC/\CC_0\ra \cohproj A$$
such that $\Phi(X)=\Ga_{\le m}(X)\mod \coh^b A$.  
To construct a functor in an opposite direction
let us first consider for every $M\in\coh A$
the functor 
$$h_M:\CC\ra \Vect_k: X\mapsto \Hom_{\GR A}(M,\Ga_*X),$$ 
where $\Vect_k$
is the category of vector spaces over $k$. We claim that $h_M$ is representable.
Indeed, for every $i\in\Z$ the functor $h_{P_i}$ is represented by $E_i$.
Now every coherent module $M$ has a presentation in the
form $M=\coker(P\ra Q),$ where $P,Q\in \PP$. We have an exact sequence of
functors $0\ra h_M\ra h_Q\ra h_P$, where the functors $h_P$ and $h_Q$
are representable. Therefore, $h_M$ is also representable.
Let us denote by $\Ga^*M$ the unique representing object for $h_M$.
Note that the functor $\Ga^*:\coh A\ra\CC$ is right exact.

For every $M\in \coh A$ there is a natural homomorphism of $A$-modules
$M\ra \Ga_{\le m}\Ga^*M$, where $M$ is concentrated in degrees $\le m$. 
We claim that the induced morphism in
$\cohproj A$ is an isomorphism. Indeed, clearly this is true for $M=P_i$. 
Now using the fact that every coherent module is a quotient of a
module from $\PP$ one can easily prove the required statement by
diagram chasing (one should check surjectivity first - see \cite{AV} (3.13)(i),(iii)).

It follows that $\Ga^*M\in\CC_0$ for $M\in \coh^b A$, so $\Ga^*$ induces a
functor
$$\Psi:\cohproj A=\coh A/ \coh^b A\ra \CC/\CC_0$$
such that $\Phi\circ \Psi\simeq \Id$. 
Furthermore, for every $X\in\CC$ and $m\in\Z$
there is a canonical morphism $\Ga^*\Ga_{\le m}X\ra X$ induced by
the embedding $\Ga_{\le m}X\ra \Ga_*X$ (considered as an element
in $\Hom(\Ga_{\le m} X,\Ga_* X)\simeq\Hom(\Ga^*\Ga_{\le m} X,X)$). 
This gives a natural
transformation of functors $\alpha:\Psi\circ \Phi \ra \Id$. We claim
that $\alpha$ is an isomorphism. Indeed, let us denote by $\PP'$ the collection
of all objects of $\CC/\CC_0$ of the form $\oplus_{j=1}^s E_{i_j}$.
One can easily check that $\a_X:\Psi(\Phi(X))\ra X$ is an isomorphism in 
$\CC/\CC_0$ for $X\in\PP'$.
Now we claim that for every $X\in\CC/\CC_0$ there exists a surjection
$f:P\ra X$ in $\CC/\CC_0$ with $P\in\PP'$. Indeed, as we have seen 
in the proof of Lemma \ref{cohdeflem}(ii), for every $X$ there exists
a morphism of the form $f:\oplus_{j=1}^s E_{i_j}\ra X$ in $\CC$ with
$\coker(f)\in\CC_0$. But such a morphism induces a surjection in $\CC/\CC_0$.
Therefore, every $X\in\CC/\CC_0$ can be represented in the form
$\coker(P\ra Q)$ with $P,Q\in \PP'$. Since $\Psi$ is right exact, this implies
that $\a_X$ is an isomorphism for all $X\in\CC/\CC_0$.
\ed

\begin{rems}
1. It should be not difficult to prove that the algebra $A(\EE)$ corresponding to an ample
sequence $\EE$ satisfies an analogue of the condition $\chi_1$ of \cite{AZ}.
Conversely, this condition for a coherent algebra $A$ should imply that 
$(P_i)$ is a projective (and hence, ample) sequence in $\cohproj A$. 
We leave the details for the reader.

2. It is easy to adapt the above theorem to the framework of graded algebras. 
Namely, if the category $\CC$ is equipped with an
autoequivalence $\si$ such that $\si(E_i)=E_{i+1}$, then the corresponding
$\Z$-algebra comes from a graded algebra. 

\noindent
3. It should be possible to weaken the assumption that all spaces 
$\Hom(E_i,X)$ are finite-dimensional in the same way as it is done in Theorem 4.5 of
\cite{AZ}. We leave for the reader to explore this. 
\end{rems}

In the remainder of the paper
we leave the general context of $\Z$-algebras and work exclusively with graded algebras. 
Theorem \ref{mainthm} can be applied to derive the following result about Veronese
subalgebras. 

\begin{prop}\label{verprop} 
Let $A=\oplus_{i\ge 0}A_i$ be a graded coherent algebra, 
$A^{(n)}=\oplus A_{in}$ be its Veronese subalgebra for some $n>0$.
Assume that $A$ is generated by $A_1$ over $k$. Then
$A^{(n)}$ is coherent and $\cohproj A\simeq \cohproj A^{(n)}$.
\end{prop}

\Pf . The idea is to apply Theorem \ref{mainthm} to the sequence $(A(in), i\in\Z)$ in
$\coh A$. First, we claim that this sequence is coherent. Indeed, it is
enough to check that for every coherent $A$-module $M$ and every
$m\in\Z$ the $A^{(n)}$-module $\oplus_{i\ge m}M_{in}$ is finitely
generated. Clearly, it is sufficient to check this for $M=A(j)$, in
which case this follows from the assumption that $A$ is generated by
$A_1$ over $k$. It remains to prove that if $M_{in}=0$ for $i>>0$,
where $M$ is a coherent $A$-module, then $M_i=0$ for $i>>0$. To this end
we observe that every element $x\in M$ of such a module satisfies
$x\cdot A_i=0$ for $i>>0$
(here we use the condition that $A$ is generated by $A_1$).
Since $M$ is finitely generated it follows that $M$ is
finite-dimensional. 
\ed

\begin{prop}\label{verprop2} 
Let $A$ be a graded algebra generated by $A_1$ over $k$ 
with a finite number of defining relations. Then $A$ is
coherent iff $A^{(n)}$ is coherent.
\end{prop}

\Pf . The ``only if" part follows from Proposition \ref{verprop}.
To prove the ``if" part by Proposition \ref{homcohprop} it suffices to 
verify that $A$ is finitely presented as a right
$A^{(n)}$-module. Note that there is a direct sum decomposition 
$$A=\oplus_{m=0}^{n-1} P^m$$ 
in the category of $A^{(n)}$-modules,
where $P^m=\oplus_i A_{m+in}$. Thus, it is enough to check that
all the $A^{(n)}$-modules $P^m$, $m=0,\ldots,n-1$, are finitely presented. 
First, we observe that they are finitely generated. Indeed, the fact that
$A$ is generated by $A_1$ immediately implies that the
$A^{(n)}$-module $P^m$ is generated by $A_m\sub P^m$. 
Since there is a finite number of defining relations between
generators of degree $1$ in $A$, we have an exact sequence of
$A$-modules
$$\oplus_{i=2}^d V_i\ot A(-i)\ra A_1\ot A(-1)\ra A\ra k\ra 0,$$
where $V_i$ are finite-dimensional vector spaces.
Therefore, we obtain exact sequences of $A^{(n)}$-modules of the form
$$\oplus_i V_i\ot P^{m-i}(-a_{m,i})\ra A_1\ot P^{m-1}\ra P^m\ra 0,$$
$m=1,\ldots, n-1$ (where $a_{m,i}\in\Z$). Now we use a simple observation that for any
homomorphism of modules $f:M_1\ra M_2$, such that $M_1$ is finitely
generated and $M_2$ is finitely presented, the module $\coker f$ is
finitely presented. Therefore, from the above sequences 
we can derive by induction in $m$ that $P^m$ is finitely presented for 
$0\le m\le n-1$. 
\ed

\begin{rem} The condition that $A$ has a finite number of defining relations
in Proposition \ref{verprop2} cannot be omitted. Indeed, let $A$ be the algebra
generated by $x$ and $y$ with the defining relations
$x^2y=0$, $yx^2=0$, $yxy=0$, $xy^{2n+1}x=0$ for $n\ge 0$. Clearly, $A$ is not coherent.
On the other hand, it is easy to
see that the algebra $A^{(2)}$ is defined by a finite number of monomial relations,
hence it is coherent by the result of \cite{Pi1}. 
\end{rem}

\section{Examples of coherent and noncoherent algebras}
\label{exsec}

Throughout this section $A$ denotes a {\it finitely generated} graded algebra.

We'll start with two reformulations of the coherence condition for such an algebra.

\begin{prop}\label{cohcrit} The following conditions are equivalent: 

\noindent
(i) $A$ is right coherent;

\noindent
(ii) for every finitely generated right (homogeneous) ideal $J\subset A$ the space
$\T^A_1(J,k)$ is finite-dimensional;

\noindent
(iii) for every $M\in \GR A$  such that the spaces $\Tor^A_0(M,k)$ and 
$\Tor^A_1(M,k)$ are finite-dimensional, the space $\Tor^A_2(M,k)$ is also finite-dimensional.
\end{prop}

\Pf . Note that the condition $\dim_k \T^A_0(M,k)<\infty$ simply means 
that $M$ is finitely generated. 

\noindent
(iii)$\implies$(ii) This follows from the isomorphism
$\T^A_i(J,k)\simeq \T^A_{i+1}(A/J,k)$ for all $i\ge 0$. 

\noindent
(ii)$\implies$(i) Since $A$ is finitely generated over $k$, it is enough to check that
for every homomorphism $f:P\ra A$, where $P\in \PP$, the module $\ker f$
is finitely generated. Let $J$ be the image of $f$. Then $J$ is finitely generated, hence
$\T^A_1(J,k)$ is finite-dimensional. Now
from the exact sequence
$$0\ra \T^A_1(J,k)\ra \T^A_0(\ker f,k)\ra \T^A_0(P,k)\ra\ldots$$
we derive that $\T^A_0(\ker f,k)$ is finite-dimensional. 
Hence, $\ker f$ is finitely generated.

\noindent
(i)$\implies$(iii) Let $M$ be a module with
finite-dimensional $\T^A_0(M,k)$ and $\T^A_1(M,k)$. Then the minimal free
resolution of $M$ has form $$\ldots\ra P^1\ra P^0\ra M\ra 0$$ where
$P^0,P^1\in \PP$. It follows that $M$ is a coherent $A$-module, so
it has resolution by modules in $\PP$ (see Lemma \ref{reslem}). Hence 
$\dim_k\T^A_i(M,k)<\infty$ for all $i$. 
\ed

\begin{cor}\label{tensorcor} 
Let $V$ be a finite-dimensional vector space over $k$.
Then the tensor algebra $T(V)$ is coherent.
\end{cor}

\Pf . Indeed, in this case $\T^A_2(M,k)=0$ for every $M$. 
\ed

It is not difficult to see that a finitely generated
$T(V)$-module is coherent if and only if its sufficiently far
truncation is a free $T(V)$-module.
This implies that every object of the category $\cohproj T(V)$ is
isomorphic to an object coming from $\PP$.
However, the morphism spaces
in $\cohproj T(V)$ between $P_i$ and $P_j$ are infinite-dimensional.
An amusing observation is that
for every $i\in\Z$ one has $P_i\simeq P_{i-1}^{\oplus\dim V}$ in
$\cohproj T(V)$.

The next two propositions provide some examples of coherent
and noncoherent algebras.

\begin{prop}\label{notprop} Assume that $A$ has a decomposition 
$A=I\oplus B$, where $I$ is a homogeneous left ideal, $B$ is a graded
subalgebra. Assume also that $B$ is a right Noetherian ring and that $I$ is
free as a left $A$-module. Then $A$ is right coherent.  
\end{prop}

\Pf . We are going to use criterion (iii) of Proposition \ref{cohcrit}. 
For every right $A$-module $M$ there is spectral sequence
with $E^2_{p,q}=\T^B_q(\T^A_p(M,B),k)$ converging to $\T^A_*(M,k)$.
On the other hand, the exact sequence $0\ra I\ra A\ra B\ra 0$ of left $A$-modules
shows that $\T^A_p(M,B)=0$ for $p\ge 2$. Now assume that $\T^A_0(M,k)$ and
$\T^A_1(M,k)$ are finite-dimensional. Then $\dim_k\T^B_0(\T^A_0(M,B),k)<\infty,$
hence, $\T^A_0(M,B)$ is a finitely generated (right) $B$-module. This implies that
all spaces $E^2_{0,q}$ are finite-dimensional. Therefore, from the assumption
$\dim_k \T^A_1(M,k)<\infty$ we derive that $E^2_{1,0}=\T^B_0(\T^A_1(M,B),k)$
is finite-dimensional. Hence, $\T^A_1(M,B)$ is a finitely generated $B$-module, which
implies that all spaces $E^2_{1,q}$ are also finite-dimensional.
\ed

\begin{cor} Let $B=k\oplus B_+$ be a graded algebra and set  
$A=B\lan z\ran/(zB_+)$, 
where $B\lan z\ran$ is the free product of $B$ with $k[z]$, $\deg(z)=1$. 
If $B$ is right Noetherian then $A$ is right coherent.
If $B$ is infinite-dimensional then $A$ is not right Noetherian.
\end{cor}

\Pf . To prove the first assertion
apply the above proposition to the decomposition $A=Az\oplus B$.
If $V_1\subset V_2\subset\ldots$ is a strictly decreasing sequence
of vector subspaces in $B_+$ then the right ideal $\oplus_{n\ge 1}V_nz^n$
in $A$ is infinitely generated.
\ed

\begin{prop}\label{noncohprop} 
Let $B=\oplus_{i\ge 0}B_i$ be a graded algebra 
equipped with a decomposition $B=I\oplus C$, where $I$ is a homogeneous right ideal and $C$ is a
graded subalgebra. Let us define the algebra $A$ as the quotient of
the algebra $B\lan z\ran$ (where $z$ has degree 1) by the relations $zI=0$ and 
$zc=cz$ for all $c\in C$. 
Assume that $I$ is not finitely generated as a right $B$-module. Then the algebra
$A$ is not right coherent.
\end{prop}

\Pf . We have a decomposition 
$A=\oplus_{n\ge 0}Bz^n$ and an exact sequence of left $A$-modules
$$0\ra A\stackrel{i}{\ra} A\ra B\ra 0,$$
where $i(x)=x\cdot z$. 
Let $J\sub A$ be the right ideal generated by $z$, so that 
$J=\oplus_{n\ge 1}Cz^n$. By Proposition \ref{cohcrit} it suffices to check
that $\dim_k\T^A_1(J,k)=\infty$.
Consider the spectral sequence with
$E^2_{p,q}=\T^B_q(\T^A_p(J,B),k)$ converging to $\T^A_*(J,k)$.
The above exact sequence shows that
$\T^A_i(J,B)\simeq \T^A_{i+1}(A/J,B)=0$ for $i\ge 1$, so
this spectral sequence degenerates in the term $E_2$ and we obtain
an isomorphism 
$$\T^A_1(J,k)\simeq \T^B_1(J\otimes_A B,k)\simeq\T^B_1(C,k).$$
But $C\simeq B/I$, so $\T^B_1(C,k)\simeq\T^B_0(I,k)$ which has
infinite dimension.
\ed

\begin{exs} 1. The algebra $A$ with generators $x,y,z$ and relations $xy=0, yz=0, xz=zx$ is
neither right nor left coherent. Indeed, this follows from Proposition \ref{noncohprop}
since we can represent $A$ in the form 
$B\lan x\ran/(xI,xz-zx)$, 
where $B=k\lan y,z\ran/(yz)$, $I=ByB=\oplus_{n\ge 0} z^nyB$, and in the form
$B'\lan z\ran/(I'z,xz-zx)$, where $B'=\lan x,y\ran/(xy)$, $I'=B'yB'=\oplus_{n\ge 0} B'yx^n$.

\noindent
2. The algebra $A$ with generators $x,y,z$ and relations
$yz=0, xz=zx$ is not left coherent (since the two-sided ideal generated by $y$ in
$k\lan x,y\ran$ is infinitely generated as a left ideal). On the other
hand, we claim that it is right coherent. 
Indeed, we can apply Proposition \ref{notprop} to the decomposition
$A=I\oplus k[x,z]$, where $I$ is the two-sided ideal generated by $y$.  
The fact that $I$ is a free left $A$-module follows from the decomposition
$I=Ay\oplus Ayx\oplus Ayx^2\oplus \ldots$. 
\end{exs}

\end{document}